\documentclass[preprint,10pt]{elsarticle}

\usepackage{amsmath, amssymb}
\usepackage{graphicx}
\journal{Applied Mathematics and Computation}

\begin{document}

\begin{frontmatter}


\title{Spatial boundary problem with the Dirichlet-Neumann condition for a singular elliptic equation}


\author{M.S.Salakhitdinov}

\address{Institute of Mathematics and Information Technologies, Uzbek Academy of Sciences, Durmon yuli str.,29, 100125 Tashkent, Uzbekistan}
\ead{salakhitdinovms@yahoo.com}
\ead[url]{http://salakhitdinov.narod.ru}

\author{E.T.Karimov}
\address{Institute of Mathematics and Information Technologies, Uzbek Academy of Sciences, Durmon yuli str.,29, 100125 Tashkent, Uzbekistan}
\ead{erkinjon@gmail.com}
\ead[url]{http://karimovet.narod.ru}

\begin{abstract}
The present work devoted to the finding explicit solution of a boundary problem with the Dirichlet-Neumann condition for elliptic equation with singular coefficients in a quarter of ball. For this aim the method of Green's function have been used. Since, found Green's function contains a hypergeometric function of Appell, we had to deal with decomposition formulas, formulas of differentiation and some adjacent relations for this hypergeometric function in order to get explicit solution of the formulated problem.

\end{abstract}


\end{frontmatter}

\section{Introduction}
\label{int}
Partial differential equations with singularities (PDEwS) have numerous applications in real life processes [1,2]. Most famous equation of this kind is Chaplygin's equation [3], which can be written as
$$
K(y)u_{xx}+u_{yy}=0
$$
or in some particular cases as
$$
u_{xx}+u_{yy}+\frac{1}{3y}u_y=0.
$$
Latter one called as Tricomi's equation and was studied well by many authors [4-7]. Due to possible applications and natural mathematical interests for generalization, (PDEwS) were studied intensively and investigations are still going on. One can find bibliographic information on this in the monographs [8-10].
Omitting huge amount works, devoted to investigations of boundary problems and potentials for two dimensional cases of aforementioned equations, note works by A.V.Bitsadze [11], A.M.Nakhushev [12], M.S.Salakhitdinov and B.Islomov [13], where three dimensional mixed type equations with singularities were investigated reducing them into two dimensional case using Fourier transformation.

Dirichlet and Dirichlet-Neumann problems for elliptic equation with one singular coefficient in some part of ball were investigated by C.Agostinelli [14], M.N.Olevskii [15]. Recently, I.T.Nazipov published a paper devoted to the investigation of the Tricomi problem in a mixed domain consisting of hemisphere and cone [16]. We also note works [17-19], where fundamental solutions and boundary problems for elliptic equations with three singular coefficient were subject of investigation. By J.J.Nieto and E.T.Karimov [20] the Dirichlet problem for an equation
$$
H_{\alpha,\beta}(u)\equiv u_{xx}+u_{yy}+u_{zz}+\frac{2\alpha}{x}u_x+\frac{2\beta}{y}u_y=0,\,\,\,0<2\alpha,2\beta<1 \eqno (1)
$$
was studied in some part of ball.

In the present work handling the method of Green's function we find explicit solution of a boundary problem with the Dirichlet-Neumann condition  for elliptic equation with two singular coefficients in a quarter of a ball.

\section{Main result}

We consider Eq.(1) in a domain
$$
\Omega=\left\{(x,y,z): x^2+y^2+z^2<R^2,\,x>0,y>0,-R<z<R\right\}
$$
which is a quarter of a ball.
\begin{figure}
\begin{center}
  \includegraphics[width=0.3\textwidth]{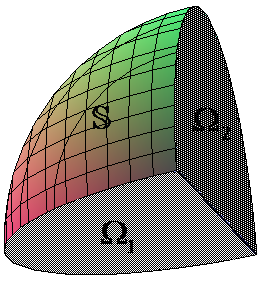}\\
  \caption{Domain $\Omega$}\label{fig1}
\end{center}
\end{figure}

\textbf{Problem D-N.} \emph{Find a function} $u\left( x,y,z \right) \in C\left( \overline \Omega  \right) \cap {C^2}\left( \Omega  \right)$,\emph{ satisfying Eq.(1) in $\Omega$ and conditions}
$$
u\left(x,y,z\right)|_{x=0} = {\tau _1}\left( y,z \right),\,\left( y,z \right) \in \overline {\Omega _1} ,\eqno (2)
$$
$$
y^{2\beta} u\left( x,y,z \right)|_{y=0} = \nu_2\left( x,z \right),\,\,\,\left( x,z \right) \in \Omega_2 ,\eqno (3)
$$
$$
u\left( x,y,z \right) = \varphi \left( x,y,z \right),\,\,\,\left( x,y,z \right) \in \overline{\mathbb{S}}.\eqno (4)
$$
Here
$$
\begin{array}{l}
\Omega_1=\left\{(x,y,z): y^2+z^2<R^2,\,x=0,y>0,-R<z<R\right\},\hfill\\
\Omega_2=\left\{(x,y,z): x^2+z^2<R^2,\,x>0,y=0,-R<z<R\right\},\hfill\\
\mathbb{S}=\left\{(x,y,z): x^2+y^2+z^2=R^2,\,x>0,y>0,-R<z<R\right\},\hfill\\
\end{array}
$$
${\tau _1}\left( y,z \right),{\nu _2}\left( x,z \right),\varphi \left( x,y,z \right)$ are given continuous functions fulfilling matching condition
$
\left. {\tau _1\left( y,z \right)} \right|_{\Upsilon _1} = \left. {\varphi \left( x,y,z \right)} \right|_{\Upsilon _1},
$ $\left(\Upsilon _1:= y^2+z^2=R^2\right)$.

\textbf{Remark.} The uniqueness of solution of the problem D-N can be proved using classical "a-b-c" method as in the work [20] or using extremume principle for elliptic equations [8]. We note also that the uniqueness theorem can be proved for more general domains.

First we give a definition of the Green's function of the problem D-N.

\textbf{Definition.} We call the function $G\left(M,M_0\right)$ as Green's function of the problem D-N, if it satisfies the following conditions:
\begin{enumerate}

\item it satisfies equation
$$
G_{xx}+G_{yy}+G_{zz}+\frac{2\alpha}{x}u_x+\frac{2\beta}{y}u_y=-\delta\left(M,M_0\right)
$$
in $\Omega$;

\item it satisfies boundary conditions
$$
\left.G\right|_{x=0}=0,\,\,y^{2\beta}\left.G_y\right|_{y=0}=0,\,\,\left.G\right|_{\overline{\mathbb{S}}}=0;
$$

\item it can be represented as
$$
G\left( M,M_0\right)=q\left( M,M_0\right)+q^*\left( M,\overline{M_0}\right). \eqno (5)
$$
\end{enumerate}
Here $\delta$ is the Dirac's delta-function, $M(x,y,z)$ is any and $M_0(x,y,z)$ is fixed point of $\Omega$ and $\overline{M_0}\left(\overline{x_0},\overline{z_0},\overline{z_0}\right)$ is a point symmetric to $M_0$ regarding the spherical surface of considered domain, i.e.
$$
\overline{x_0}=-\frac{R^2}{R_0^2}x_0,\,\,\overline{y_0}=-\frac{R^2}{R_0^2}y_0,\,\,\overline{z_0}=-\frac{R^2}{R_0^2}z_0,\,\,R_0^2=x_0^2+y_0^2+z_0^2,
$$
$$
q\left(M,M_0\right)=k\left(r^2\right)^{\alpha-\beta-\frac{3}{2}}\left(xx_0\right)^{1-2\alpha}F_2\left(\frac{3}{2}-\alpha+\beta,1-\alpha,
\beta;2-2\alpha,2\beta;\xi,\eta\right) \eqno (6)
$$
is one of the fundamental solutions of Eq.(1) [20],
$$
{F_2}\left( {a;{b_1},{b_2};{c_1},{c_2};x,y} \right) = \sum\limits_{i,j = 0}^\infty  {\frac{{{{\left( a \right)}_{i + j}}{{\left( {{b_1}} \right)}_i}{{\left( {{b_2}} \right)}_j}}}{{{{\left( {{c_1}} \right)}_i}{{\left( {{c_2}} \right)}_j}i!j!}}{x^i}{y^j}}
$$
is Appel's hypergeometric function [21],
$$
q^*\left(M,M_0\right)=-\left(\frac{R}{R_0}\right)^{3-2\alpha+2\beta}q\left(M,\overline{M}_0\right)
$$
is a regular part of the Green's function (5), i.e. satisfies Eq.(1) in any point of $\Omega$,
$$
\xi=1-\frac{r_1^2}{r^2},\,\,\eta=1-\frac{r_2^2}{r^2},\,\,r^2=(x-x_0)^2+(y-y_0)^2+(z-z_0)^2,
$$
$$
r_1^2=(x+x_0)^2+(y-y_0)^2+(z-z_0)^2,\,\,r_2^2=(x-x_0)^2+(y+y_0)^2+(z-z_0)^2,
$$
$$
k=\frac{1}{2\pi}\frac{\Gamma(1-\alpha)\Gamma(\beta)\Gamma(2-2\alpha+2\beta)}{\Gamma(2-2\alpha)\Gamma(2\beta)\Gamma(1-\alpha+\beta)}.
\eqno (7)
$$

Excise from the domain $\Omega$ a small ball with its center at $M_0$ and with radius $\rho>0$, designating the sphere of the excised ball as $C_\rho$ and by $\Omega_\rho$ denote the remaining part of $\Omega$.

Applying Green's formula [22] in this domain, one can obtain the following:
$$
\begin{array}{l}
\iint\limits_{C_\rho}{x^{2\alpha}y^{2\beta}\left[\displaystyle{u\frac{\partial G}{\partial n}-G\frac{\partial u}{\partial n}}\right]ds}=\iint\limits_{\Omega_1}y^{2\beta}\tau_1(y,z)G^*\left(M,M_0\right)dydz\hfill\\
+\iint\limits_{\Omega_2}x^{2\alpha}\nu_2(x,z)G^{**}\left(M,M_0\right)dxdz+
\iint\limits_{\mathbb{S}}x^{2\alpha}y^{2\beta}\varphi(\sigma)\displaystyle{\frac{\partial G}{\partial n}}d\sigma.\hfill\\
\end{array}
\eqno (8)
$$
Here $G^*\left(M,M_0\right)=\left.x^{2\beta}\displaystyle{\frac{\partial G\left(M,M_0\right)}{\partial x}}\right|_{x=0}$, $G^{**}\left(M,M_0\right)=\left.G\left(M,M_0\right)\right|_{y=0}$, $n$ is outer normal to the $\partial \Omega$.

Let us first evaluate an integral
$$
I=\iint\limits_{C_\rho}x^{2\alpha}y^{2\beta}u\frac{\partial G}{\partial n}ds. \eqno (9)
$$
Based on (5) we write (9) as follows:
$$
I=I_1+I_2=\iint\limits_{C_\rho}x^{2\alpha}y^{2\beta}u\frac{\partial q}{\partial n}ds+\iint\limits_{C_\rho}x^{2\alpha}y^{2\beta}u\frac{\partial q^*}{\partial n}ds,
$$
where
$$
\frac{\partial q}{\partial n}=\frac{\partial q}{\partial x}\cos(n,x)+\frac{\partial q}{\partial y}\cos(n,y)+\frac{\partial q}{\partial z}\cos(n,z). \eqno (10)
$$
Formal calculation gives
$$
\begin{array}{l}
\displaystyle{\frac{\partial q}{\partial x}=\frac{\partial}{\partial x}\left[kP_1P_2F_2(...)\right]}\hfill\\
\displaystyle{=k\left[\left(\frac{\partial P_1}{\partial x}P_2+P_1\frac{\partial P_2}{\partial x}\right)F_2(...)+P_1P_2\left(\frac{\partial F_2}{\partial \xi}\xi_x+\frac{\partial F_2}{\partial \eta}\eta_x\right)\right]},\hfill\\
\end{array}
$$
where $P_1=\left(r^2\right)^{\alpha-\beta-\frac{3}{2}},\,\,P_2=\left(xx_0\right)^{1-2\alpha}$.

Using formula of differentiation for hypergeometric function $F_2(...)$ [21]
$$
\begin{array}{l}
  \displaystyle{\frac{\partial ^{i + j}F_2\left(a,b_1,b_2;c_1,c_2;x,y \right)}{\partial x^iy^j} }=  \hfill \\
  \displaystyle{\frac{\left( a \right)_{i + j}\left(b_1\right)_i\left(b_2\right)_j}{\left(c_1\right)_i\left(c_2\right)_j}F_2\left(a + i + j,b_1 + i,b_2 + j;c_1 + i,c_2 + j;x,y \right)} \hfill \\
\end{array}
$$
we get
$$
\begin{array}{l}
  \displaystyle{\frac{\partial}{\partial\xi}F_2\left(\frac{3}{2}-\alpha+\beta,1-\alpha,\beta;2-2\alpha,2\beta;\xi,\eta\right)}\hfill \\
 \displaystyle{= \frac{\left(\frac{3}{2}-\alpha+\beta \right)\left(1-\alpha \right)}{\left(2-2\alpha\right)}F_2\left(\frac{5}{2}-\alpha+\beta,2-\alpha,\beta;3-2\alpha,2\beta;\xi,\eta \right)}, \hfill \\
\end{array}
$$
$$
\begin{array}{l}
  \displaystyle{\frac{\partial}{\partial\eta}F_2\left(\frac{3}{2}-\alpha+\beta,1-\alpha,\beta;2-2\alpha,2\beta;\xi,\eta\right)}\hfill \\
 \displaystyle{= \frac{\left(\frac{3}{2}-\alpha+\beta \right)\beta}{2\beta}F_2\left(\frac{5}{2}-\alpha+\beta,1-\alpha,\beta+1;2-2\alpha,2\beta+1;\xi,\eta \right).} \hfill \\
\end{array}
$$
Considering
$$
\frac{\partial P_1}{\partial x}=\frac{2\left(x-x_0\right)}{r^2}\left(\alpha-\beta-\frac{3}{2}\right)P_1,\,\,\frac{\partial P_2}{\partial x}=\frac{1-2\alpha}{x}P_2,
$$
$$
\xi_x=-\frac{4x_0}{r^2}-\frac{2\left( x-x_0\right)}{r^2}\xi,\,\eta_x=-\frac{2\left( x-x_0\right)}{r^2}\eta,
$$
we have
$$
\begin{array}{l}
\displaystyle{\frac{\partial q}{\partial x}=\frac{kP_1P_2}{r^2}\left[2\left(x-x_0\right)\left(\alpha-\beta-\frac{3}{2}\right)\right.}\hfill\\
\displaystyle{\times F_2\left(\frac{3}{2}-\alpha+\beta,1-\alpha,\beta;2-2\alpha,2\beta;\xi,\eta\right)+\frac{1-2\alpha}{x}r^2}\hfill\\
\displaystyle{\times F_2\left(\frac{3}{2}-\alpha+\beta,1-\alpha,\beta;2-2\alpha,2\beta;\xi,\eta\right)-2x_0\left(\frac{3}{2}-\alpha+\beta\right)}\hfill\\
\displaystyle{\times F_2\left(\frac{5}{2}-\alpha+\beta,2-\alpha,\beta;3-2\alpha,2\beta;\xi,\eta\right)-2\left(x-x_0\right)\left(\frac{3}{2}-\alpha+\beta\right)
\xi\frac{1-\alpha}{2-2\alpha}}\hfill\\
\displaystyle{\times F_2\left(\frac{5}{2}-\alpha+\beta,2-\alpha,\beta;3-2\alpha,2\beta;\xi,\eta\right)-2\left(x-x_0\right)\left(\frac{3}{2}-\alpha+\beta\right)
\eta\frac{\beta}{2\beta}}\hfill\\
\displaystyle{\left.\times F_2\left(\frac{5}{2}-\alpha+\beta,1-\alpha,\beta+1;2-2\alpha,2\beta+1;\xi,\eta\right)\right]}.\hfill\\
\end{array}
\eqno (11)
$$
Now we apply to (11) the following adjacent relation [21]
$$
\begin{array}{l}
  \displaystyle{x\frac{b_1}{c_1}F_2\left(a + 1,b_1 + 1,b_2;c_1 + 1,c_2;x,y \right)} \hfill \\
  \displaystyle{  + y\frac{b_2}{c_2}F_2\left(a + 1,b_1,b_2 + 1;c_1,c_2 + 1;x,y \right)} \hfill \\
    \displaystyle{ = F_2\left( a + 1,b_1,b_2;c_1,c_2;x,y \right) - F_2\left(a,b_1,b_2;c_1,c_2;x,y \right) }\hfill \\
\end{array}
$$
and will obtain
$$
\begin{array}{l}
  \displaystyle{\frac{\partial q}{\partial x}=\frac{kP_1P_2}{r^2}\left[2\left(x-x_0\right)\left(\alpha-\beta-\frac{3}{2}\right)\right.}\hfill\\
    \displaystyle{\times F_2\left(\frac{5}{2}-\alpha+\beta,1-\alpha,\beta;2-2\alpha,
  2\beta;\xi,\eta\right)}\hfill\\
  \displaystyle{\left.+2x_0\left(\alpha-\beta-\frac{3}{2}\right)F_2\left(\frac{5}{2}-\alpha+\beta,2-\alpha,\beta;3-2\alpha,2\beta;\xi,\eta\right)
  \right]}\hfill\\
  \displaystyle{+kP_1P_2\frac{1-2\alpha}{x}r^2F_2\left(\frac{3}{2}-\alpha+\beta,1-\alpha,\beta;2-2\alpha,2\beta;\xi,\eta\right)}.\hfill\\
\end{array}
\eqno (12)
$$
Similarly we get
$$
\begin{array}{l}
  \displaystyle{\frac{\partial q}{\partial y}=\frac{kP_1P_2}{r^2}\left[2\left(y-y_0\right)\left(\alpha-\beta-\frac{3}{2}\right)\right.}\hfill\\
  \displaystyle{ \times F_2\left(\frac{5}{2}-\alpha+\beta,1-\alpha,\beta;2-2\alpha,2\beta;\xi,\eta\right)}\hfill\\
  \displaystyle{\left.+2y_0\left(\alpha-\beta-\frac{3}{2}\right)F_2\left(\frac{5}{2}-\alpha+\beta,1-\alpha,\beta+1;2-2\alpha,2\beta+1;\xi,
  \eta\right)\right]},\hfill\\
\end{array}
\eqno (13)
$$
$$
\frac{\partial q}{\partial z}=\frac{kP_1P_2}{r^2}2\left(z-z_0\right)\left(\alpha-\beta-\frac{3}{2}\right)F_2\left(\frac{5}{2}-\alpha+\beta,1-\alpha,\beta;2-2\alpha,
  2\beta;\xi,\eta\right).
\eqno (14)
$$
Taking (12)-(14) into account, from (10) we get the following:
$$
\begin{array}{l}
    \displaystyle{\frac{\partial }{\partial n}q\left( M,M_0 \right)
   = \left( \alpha  - \beta  - \frac{3}{2} \right)k P_1 P_2 }\hfill\\
     \displaystyle{\times F_2\left(\frac{5}{2} - \alpha  + \beta,1 - \alpha ,\beta ;2 - 2\alpha, 2\beta;\xi,\eta\right)\frac{\partial }{\partial n}\left[ \ln {r^2} \right]+ \frac{k P_1 P_2}{r^2}\left( \alpha  - \beta  - \frac{3}{2} \right)}\hfill\\
       \displaystyle{\times\left[ 2x_0F_2\left(\frac{5}{2} - \alpha  + \beta,2 - \alpha ,\beta ;3 - 2\alpha ,2\beta ;\xi ,\eta \right)\right.} \hfill \\
    \displaystyle{\left. + 2y_0F_2\left( \frac{5}{2} - \alpha  + \beta,1 - \alpha ,\beta+1 ;2 - 2\alpha ,2\beta+1 ;\xi ,\eta \right)\right]}\hfill\\
       \displaystyle{+ k P_1 P_2\frac{1 - 2\alpha}{x}F_2\left(\frac{3}{2} - \alpha  + \beta,1 - \alpha ,\beta ;2 - 2\alpha ,2\beta ;\xi ,\eta \right).} \hfill \\
\end{array}
\eqno (15)
$$

Now consider the integral
$$
\begin{array}{l}
  \displaystyle{I_1 = \iint\limits_{C_\rho} {x^{2\alpha }y^{2\beta}u\frac{\partial q}{\partial n}dS} = I_{11}+I_{12}+I_{13}}\hfill\\
  \displaystyle{=\iint\limits_{C_\rho}{x^{2\alpha }y^{2\beta}\left( \alpha  - \beta  - \frac{3}{2} \right)k P_1 P_2 }}\hfill\\
    \displaystyle{\times F_2\left(\frac{5}{2} - \alpha  + \beta,1 - \alpha ,\beta ;2 - 2\alpha, 2\beta;\xi,\eta\right)\frac{\partial }{\partial n}\left[ \ln {r^2} \right]u\left(M_0\right)ds}\hfill\\
  \displaystyle{+\iint\limits_{C_\rho}{x^{2\alpha }y^{2\beta}\frac{k P_1 P_2}r^2\left( \alpha  - \beta  - \frac{3}{2} \right)}}\hfill\\
    \displaystyle{\times\left[ 2x_0F_2\left(\frac{5}{2} - \alpha  + \beta,2 - \alpha ,\beta ;3 - 2\alpha ,2\beta ;\xi ,\eta \right)\right.}\hfill\\
  \displaystyle{\left. + 2y_0F_2\left( \frac{5}{2} - \alpha  + \beta,1 - \alpha ,\beta+1 ;2 - 2\alpha ,2\beta+1 ;\xi ,\eta \right)\right]u\left(M_0\right)ds}\hfill\\
  \displaystyle{+ \iint\limits_{C_\rho}{x^{2\alpha }y^{2\beta}k P_1 P_2\frac{1 - 2\alpha}{x}F_2\left(\frac{3}{2} - \alpha  + \beta,1 - \alpha ,\beta ;2 - 2\alpha ,2\beta ;\xi ,\eta \right)u\left(M_0\right)ds}}. \hfill \\
\end{array}
$$

We use the following spherical system of coordinates [22]:
$$
\begin{array}{l}
x = {x_0} + \rho \sin \theta \cos \psi ,\,\,y = {y_0} + \rho \sin \theta \sin \psi ,\,\,z = {z_0} + \rho \cos \theta ,\hfill\\
0 < \theta  < \pi ,\,0 < \psi  < 2\pi ,\,0 < \rho  < \infty .\hfill\\
\end{array}
$$
Then we have
$$
\begin{array}{l}
  I_{11} = -2\left( \alpha  - \beta  - \frac{3}{2} \right) k\int\limits_0^{2\pi } d\psi \int\limits_0^\pi  {x_0^{1 - 2\alpha }\left(x_0 + \rho \sin \theta \cos \psi \right)\left(y_0 + \rho \sin \theta \sin \psi \right)^{2\beta}}  \hfill \\
   \times u\left( x_0 + \rho \sin \theta \cos \psi ,\,y_0 + \rho \sin \theta \sin \psi ,\,z_0 + \rho \cos \theta \right)\left( \rho ^2\right)^{\alpha  -\beta  - 1}\hfill\\
  \times F_2\left(\frac{5}{2}-\alpha+\beta, 1-\alpha,\beta;2-2\alpha,2\beta;\xi_s,\eta_s\right)\sin \theta d\theta . \hfill \\
\end{array}
$$
Let us first evaluate hypergeometric function $F_2$. We use decomposition formula [21]
$$
\begin{array}{l}
  F_2\left(a,b_1,b_2;c_1,c_2;x,y \right) =  \hfill \\
  \sum\limits_{i = 0}^\infty  \frac{\left( a \right)_i\left(b_1\right)_i\left(b_2\right)_i}{\left(c_1\right)_i\left(c_2\right)_i i!}x^iy^i {}_2F_1\left(a + i,b_1 + i;c_1 + i;x \right)\cdot {}_2F_1\left( a + i,b_2 + i;c_2 + i;y \right), \hfill \\
\end{array}
$$
after then using auto transformation formula [23]
$$
_2{F_1}\left( {a,b,c;x} \right) = \left(1 - x\right)^{-b}{}_2F_1\left(c - a,b,c;\frac{x}{x - 1}\right)
$$
we obtain
$$
\begin{array}{l}
   \displaystyle{ F_2\left(\frac{5}{2}-\alpha+\beta, 1-\alpha,\beta;2-2\alpha,2\beta;\xi_s,\eta_s\right) =}\hfill\\
      \displaystyle{\times\sum\limits_{i = 0}^\infty  {\frac{\left(\frac{5}{2} - \alpha  + \beta \right)_i\left(1 - \alpha\right)_i\left(\beta\right)_i}{\left( 2 - 2\alpha \right)_i\left(2\beta\right)_i i!}\xi _s^i\eta _s^i}}\hfill\\
    \displaystyle{\times {}_2F_1\left(\frac{5}{2} - \alpha  + \beta + i ,1 - \alpha  + i;2 - 2\alpha  + i;\xi_s\right)}\hfill\\
   \displaystyle{ \times {}_2F_1\left(\frac{5}{2} - \alpha  + \beta + i ,\beta  + i; 2\beta + i;\eta_s\right)}\hfill\\
  \displaystyle{  =\sum\limits_{i = 0}^\infty  {\frac{\left(\frac{5}{2} - \alpha  + \beta \right)_i\left(1 - \alpha\right)_i\left(\beta\right)_i}{\left( 2 - 2\alpha \right)_i\left(2\beta\right)_i i!}\xi _s^i\eta _s^i}\left(1 - \xi_s\right)^{\alpha - 1 - i}\left( 1 - \eta_s\right)^{- \beta - i} } \hfill \\
  \displaystyle{  \times {}_2F_1\left(- \frac{1}{2} - \alpha - \beta ,1 - \alpha  + i;2 - 2\alpha  + i;\frac{\xi_s}{\xi_s - 1}\right)}\hfill\\
  \displaystyle{  \times {}_2F_1\left(- \frac{5}{2} + \alpha  + \beta ,\beta  + i;2\beta  + i;\frac{\eta_s}{\eta_s - 1} \right).} \hfill \\
\end{array}
$$
Here
$$
\begin{array}{l}
    \displaystyle{\frac{\xi _s}{\xi_s - 1} = 1 - \frac{\rho^2}{r_{1s}^2},\,\frac{\eta_s}{\eta_s - 1} = 1 - \frac{\rho ^2}{r_{2s}^2},\,r_{1s}^2 = 4x_0^2 + 4{x_0}\rho \sin \theta \cos \psi ,} \hfill \\
  \displaystyle{  r_{2s}^2 = 4y_0^2 + 4{y_0}\rho \sin \theta \sin \psi ,\,1 - \xi_s = \frac{r_{1s}^2}{\rho^2},\,1 - \eta_s = \frac{r_{2s}^2}{\rho ^2}.} \hfill \\
\end{array}
$$
After some evaluations we find
$$
\begin{array}{l}
\displaystyle{F_2\left(\frac{5}{2}-\alpha+\beta, 1-\alpha,\beta;2-2\alpha,2\beta;\xi_s,\eta_s\right) = \left( r_{1s}^2\right)^{\alpha  - 1}\left( r_{2s}^2\right)^{-\beta}\left(\rho^2\right)^{1 - \alpha + \beta } }\hfill \\
 \displaystyle{  \times \sum\limits_{i = 0}^\infty  {\frac{\left(\frac{5}{2} - \alpha + \beta\right)_i\left(1 - \alpha\right)_i\left(\beta\right)_i}{\left(2 - 2\alpha\right)_i\left(2\beta\right)_i i!}\left(\frac{\rho^2}{r_{1s}^2} - 1\right)^i\left(\frac{\rho^2}{r_{2s}^2} - 1\right)^i} } \hfill \\
\displaystyle{  \times {}_2F_1\left(- \frac{1}{2} - \alpha - \beta ,1 - \alpha  + i;2 - 2\alpha  + i;\frac{\xi_s}{\xi_s - 1}\right)}\hfill\\
\displaystyle{  \times {}_2F_1\left(- \frac{5}{2} + \alpha  + \beta ,\beta  + i;2\beta  + i;\frac{\eta_s}{\eta_s - 1} \right).} \hfill \\
\end{array}
$$
Considering [23]
$$
\begin{array}{l}
\displaystyle{_2F_1\left(a,b,c;1\right) = \sum\limits_{n = 0}^\infty  {\frac{\left( a \right)_n\left( b \right)_n}{\left( c \right)_n n!}}  = \frac{\Gamma \left( c \right)\Gamma \left( c - a - b \right)}{\Gamma \left( c - a \right)\Gamma \left( c - b \right)}},\hfill\\
c \ne 0, - 1, - 2,...,\,\Re \left( {c - a - b} \right) > 0.\hfill\\
\end{array}
$$
we deduce
$$
\begin{array}{l}
\displaystyle{\mathop {\lim }\limits_{\rho  \to 0} {}_2F_1\left(- \frac{1}{2} - \alpha  - \beta ,1 - \alpha  + i;2 - 2\alpha  + i;1 - \frac{\rho^2}{r_{1s}^2}\right)}\hfill\\
\displaystyle{ = \frac{\Gamma \left( 2 - 2\alpha  + i\right)\Gamma \left(\beta  + \frac{3}{2}\right)}{\Gamma \left(\frac{5}{2} - \alpha + \beta  + i \right)\Gamma \left( 1 - \alpha \right)}},\hfill\\
 \end{array}
$$
$$
\begin{array}{l}
\displaystyle{\mathop {\lim }\limits_{\rho  \to 0} {}_2F_1\left(- \frac{5}{2} + \alpha  + \beta ,\beta  + i;2\beta  + i;1 - \frac{\rho^2}{r_{2s}^2}\right)} \hfill\\
\displaystyle{= \frac{\Gamma \left( 2\beta  + i\right)\Gamma \left(- \alpha  + \frac{5}{2}\right)}{\Gamma \left( \frac{5}{2} - \alpha  + \beta  + i\right)\Gamma \left( \beta\right)}.}\hfill\\
\end{array}
$$
Now we get
$$
\begin{array}{l}
\displaystyle{\mathop {\lim }\limits_{\rho  \to 0} I_{11} =  - 2\pi  \cdot k \cdot \left( \alpha  - \beta  - \frac{3}{2}\right) \cdot 2^{2\alpha  - 2\beta+1}u\left(M_0\right)}\hfill\\
 \displaystyle{\times \frac{\Gamma \left( \frac{5}{2} - \alpha \right)\Gamma \left(\frac{3}{2} + \beta\right)}{\Gamma \left(1 - \alpha\right)\Gamma \left(\beta\right)} \cdot \mathfrak{P},}\hfill\\
 \end{array}
$$
where
$$
\mathfrak{P} = \sum\limits_{i = 0}^\infty \frac{\left(\frac{5}{2} - \alpha  + \beta \right)_i\left(1 - \alpha\right)_i\left(\beta\right)_i}{\left( 2 - 2\alpha \right)_i\left(2\beta\right)_i i!} \cdot \frac{\Gamma \left( 2 - 2\alpha  + i\right)\Gamma \left(2\beta  + i \right)}{\Gamma^2\left( \frac{5}{2} - \alpha  + \beta  + i \right)}.
$$
Considering [23]
$$
\Gamma\left(2-2\alpha+i\right)=\Gamma\left(2-2\alpha\right)\left(2-2\alpha\right)_i,\,
\Gamma\left(2\beta+i\right)=\Gamma\left(2\beta\right)\left(2\beta\right)_i,\,\Gamma\left(\frac{3}{2}\right)=\frac{\sqrt{\pi}}{2},
$$
$$
\Gamma\left(\frac{5}{2}-\alpha+\beta\right)=\frac{\sqrt{\pi}\Gamma\left(2-2\alpha+2\beta\right)\left(\frac{3}{2}-\alpha+\beta\right)}{2^{1-2\alpha+2\beta}
\Gamma\left(1-\alpha+\beta\right)}\left(\frac{5}{2}-\alpha+\beta\right)_i,
$$
making some calculations we obtain
$$
\mathfrak{P} = \frac{2^{- 2\alpha  + 2\beta}}{\left( \frac{3}{2} - \alpha  + \beta \right)} \cdot \frac{\Gamma \left( 2 - 2\alpha \right)\Gamma \left( 2\beta\right)\Gamma \left(1 - \alpha  + \beta \right)}{\Gamma \left(2 - 2\alpha  + 2\beta\right)\Gamma \left(\frac{5}{2} - \alpha\right)\Gamma \left( \frac{3}{2} + \beta\right)}.
$$
Therefore
$$
\mathop {\lim }\limits_{\rho  \to 0} I_{11} = k \cdot 2\pi  \cdot \frac{\Gamma \left(2 - 2\alpha\right)\Gamma \left(2\beta\right)\Gamma \left( 1 - \alpha  + \beta\right)}{\Gamma \left( 2 - 2\alpha  + 2\beta \right)\Gamma\left(1-\alpha\right)\Gamma\left(\beta\right)} \cdot u\left(M_0\right).
$$
If we choose $k$ as in (7), we will have
$$
\mathop {\lim }\limits_{\rho  \to 0} I_{11} = u\left(M_0\right).
$$
By similar evaluations one can get that
$$
\mathop {\lim }\limits_{\rho  \to 0} I_{12} = \mathop {\lim }\limits_{\rho  \to 0} I_{13} = \mathop {\lim }\limits_{\rho  \to 0} I_2 = 0.
$$

If we consider an integral
$$
\iint\limits_{C_\rho} {x^{2\alpha}y^{2\beta}G\frac{\partial u}{\partial n}dS},
$$
using above given algorithm for evaluations, we can prove that
$$
\mathop {\lim }\limits_{\rho  \to 0} \iint\limits_{C_\rho} {x^{2\alpha}y^{2\beta}G\frac{\partial u}{\partial n}dS = 0}.
$$

Now we can formulate our result as the following

\textbf{Theorem.} The Problem D-N has unique solution represented as follows
$$
\begin{array}{l}
  u\left(M_0\right) = \iint\limits_{\Omega _1} {y^{2\beta} \tau_1\left(y,z\right)G^*\left(M,M_0\right)dydz}
   + \iint\limits_{\Omega _2} {x^{2\alpha} \nu_2\left(x,z\right)G^{**}\left(M,M_0\right)dxdz} \hfill\\
   + \displaystyle{\iint\limits_\mathbb{S} {x^{2\alpha} y^{2\beta}\varphi \left( \sigma  \right)\frac{\partial G\left(M,M_0\right)}{\partial n}d\sigma}.} \hfill \\
\end{array}
$$

The particular values of the Green's function is given by
$$
\begin{array}{l}
  G^*\left( M,M_0\right) = k\left( 1 - 2\alpha \right)x_0^{1 - 2\alpha }\hfill \\\left[ \frac{{}_2F_1\left(\frac{3}{2} - \alpha  + \beta ,\beta ,2\beta ;\eta _{ox} \right)}{\left(r^2|_{x=0}\right)^{-\frac{3}{2} + \alpha  - \beta}}   - \left(\frac{R}{R_0} \right)^{2 - 4\alpha}\frac{{}_2F_1\left( \frac{3}{2} - \alpha  + \beta , \beta ,2\beta ;\overline {\eta}_{ox} \right)}{\left(\overline{r}^2|_{x=0}\right)^{-\frac{3}{2} + \alpha  - \beta}}\right], \hfill \\
\end{array}
$$
$$
\begin{array}{l}
G^{**}\left(M,M_0\right) = k\left(xx_0\right)^{1 - 2\alpha}\hfill\\
  \left[ \frac{{}_2F_1\left(\frac{3}{2} - \alpha  + \beta ,1 - \alpha ,2 - 2\alpha ;\xi _{oy}\right)}{\left(r^2|_{y=0}\right)^{-\frac{3}{2} + \alpha  - \beta }} - \left(\frac{R}{R_0} \right)^{2 - 4\alpha}\frac{{}_2F_1\left(\frac{3}{2} - \alpha  + \beta ,1 - \alpha ,2 - 2\alpha ;\overline {\xi}_{oy}\right)}{\left(\overline{r}^2|_{y=0}\right)^{-\frac{3}{2} + \alpha  - \beta}} \right]. \hfill \\
\end{array}
$$
Here,
$$
\begin{array}{l}
\overline{r}^2|_{x=0}=\left( R - \frac{yy_0}{R}\right)^2 + \left( R - \frac{zz_0}{R} \right)^2 + \frac{x_0^2 + z_0^2}{R^2}y^2 + \frac{x_0^2 + y_0^2}{R^2}z^2 - R^2,\hfill\\
\overline{r}^2|_{y=0}=\left( R - \frac{xx_0}{R}\right)^2 + \left( R - \frac{zz_0}{R} \right)^2 + \frac{y_0^2 + z_0^2}{R^2}x^2 + \frac{x_0^2 + y_0^2}{R^2}z^2 - R^2,\hfill\\
\xi_{0y}=-\frac{4xx_0}{\left.r^2\right|_{y=0}},\,\,\eta_{0x}=-\frac{4yy_0}{r^2|_{x=0}},\,\, \overline{\xi}_{0y}=-\frac{4R^2}{R_0^2}\frac{xx_0}{\overline{r}^2|_{y=0}}, \,\, \overline{\eta}_{0x}=-\frac{4R^2}{R_0^2}\frac{yy_0}{\overline{r}^2|_{x=0}}, \hfill \\
\end{array}
$$

\section{Acknowledgment}
Second author (Karimov E.T.) thanks ICTP (Abdus Salam International Center of Theoretical Physics) for given possibility to get this result via offering visiting fellowship.

\bibliographystyle{elsarticle-num}

\end{document}